\documentclass[preprint,authoryear,12pt]{elsarticle}
\usepackage{a4wide}
\usepackage{amsmath,amssymb,amsthm}

\newtheorem{thm}{Theorem}
\newtheorem{coro}{Corollary}
\newtheorem{assh}{Hypothesis}
\newdefinition{defi}{Definition}
\newenvironment{mypr}
{
\par\vspace{.5\baselineskip}\noindent\textbf{Proof.}}
{\nobreak\hfill$\Box$\par\vspace{.5\baselineskip}}

\usepackage{lebListes}  
\noitemsep  

\usepackage{ifthen}

\usepackage{algorithm,algpseudocode}

\usepackage{color,pstcol,pstricks}
\newgray{backgroundgray}{0.90}
\definecolor{backgroundgray}{gray}{0.90}

\newenvironment{petit}
{\par\vspace{.5\baselineskip}\noindent\footnotesize}
{\nobreak\par\vspace{.5\baselineskip}}

\def\tempfolder{./temp/}

\providecommand{\imprimeTexteCache}{non}
\providecommand{\imprimeTexteSecret}{non}

\newif\ifsol
\newif\ifnfs

\ifthenelse{\equal{\imprimeTexteCache}{oui}}{\soltrue}{\solfalse }
\ifthenelse{\equal{\imprimeTexteSecret}{oui}}{\nfstrue}{\nfsfalse}

\ifsol
   \providecommand{\sol}[1]{\bs #1 \es}
\else
   \providecommand{\sol}[1]{}
\fi

\ifnfs
   \providecommand{\nfs}[1]{~\\{\footnotesize \sffamily \emph{\textbf{Note: } #1}~\\}}
\else
   \providecommand{\nfs}[1]{}
\fi

\newenvironment{mySolution}
{\red \it
\par\vspace{.5\baselineskip}\noindent\textbf{Solution.~}}
{\vspace{.5\baselineskip}}
\def\bs{\begin{mySolution}}
\def\es{\end{mySolution}}

\def\problemfile#1{%
\begin{filecontents}{\tempfolder#1}}

\ifnfs

\else

\fi

\ifnfs

\else

\fi

\ifnfs

\else

\fi

\ifnfs
\newcommand{\tbd}[1]{\textsc{To be done:\\#1}}
\else
\newcommand{\tbd}[1]{}
\fi

\newtheorem{maquestion}{\sc Question}[section]

\ifsol

\else

\fi

\newcommand{\imp}[1]{{\blue \bf #1}}

\newcommand{\bp}{ 
  \small \ttfamily
 \begin{tabbing}
 aaa\=aaa\=aaa\=aaa\=aaaaaaaa \= aaaaaaaaaa\= \kill
 }
\newcommand{\ep}{\end{tabbing}\normalfont\normalsize }

\newcommand{\cref}[1]{Chapter~\ref{#1}}

\newcommand{\eref}[1]{Eq.(\ref{#1})}

\newcommand{\calC}{\mathcal C}

\newcommand{\calF}{\mathcal F}

\newcommand{\calP}{\mathcal P}

\newcommand{\lp}{\left(}
\newcommand{\lb}{\left[}

\newcommand{\rp}{\right)}
\newcommand{\rb}{\right]}

\newcommand{\eps}{\epsilon}
\newcommand{\esp}[1]{\E\left(#1\right)}
\newcommand{\espc}[2]{\E\left(\left.#1 \right| #2\right)}

 \newcommand{\norm}[1]{\left\|#1\right\|}

\def\be{\begin{equation}}
\def\ee{\end{equation}}
\def\ben{\[}
\def\een{\]}
\def\bearn{\begin{eqnarray*}}
\def\eearn{\end{eqnarray*}}
\def\bear{\begin{eqnarray}}
\def\eear{\end{eqnarray}}
\def\barr{\begin{array}}
\def\earr{\end{array}}

\def\bmat{\left(\begin{array}}
\def\emat{\end{array}\right)}

\newcommand{\limit}[2]{\lim_{#1 \rightarrow #2}}
\newcommand{\eqdef}{\stackrel{\mathrm{def}}{=}} 

\def\Reals{\mathbb{R}}

\def\Nats{\mathbb{N}}
\def\E{\mathbb{E}}
\def\P{\mathbb{P}}
\def\T{\mathbb{T}}

\newcommand{\bracket}[1]{ \left\{\begin{array}{l} #1 \end{array} \right.}

\newsavebox{\coloredbox}

\newsavebox{\traitbox}

\title{On Mean Field Convergence and Stationary Regime}

\author[unine]{Michel Bena\"im}
\ead{michel.benaim@unine.ch}

\author[epfl]{Jean-Yves Le Boudec}
\ead{jean-yves.leboudec@epfl.ch}
\address[unine]{Institut de Math\'ematiques, University of Neuch\^atel, Switzerland}
\address[epfl]{EPFL IC-LCA2 - Lausanne, Switzerland}

\begin{document}

\begin{frontmatter}
\begin{abstract}
  Assume that a family of stochastic processes on some Polish
space $E$ converges to a deterministic process; the convergence
is in distribution (hence in probability) at every fixed point in
time. This assumption holds for a large family of processes,
among which many mean field interaction models and is weaker than previously assumed.
We show that
any limit point of an invariant probability of the stochastic
process is an invariant probability of the deterministic
process. The results are valid in discrete and in continuous
time.

\end{abstract}
\end{frontmatter}

\section{Introduction}
%

This paper is motivated by results on mean field interaction
models and stochastic approximation algorithms which obtain
convergence of a family of stochastic  processes to a
deterministic
limit \citep{kurtz1970solutions,sandholm2006population,prout07,graham1994chaos,
benaim2008cmf}. Often, convergence is over finite time
horizons, which asks the question of whether the convergence
extends to the stationary regime. In this paper, we show that some form of
convergence of the stationary regimes follows systematically
from convergence at any fixed and finite time horizon, under
very weak assumptions. Previous answers to this question exist with stronger assumptions than here;
for example, the space is finite dimensional and the deterministic limit is a dissipative ODE \citep{benaim98}, or the set of invariant distributions is tight
\citep{fort1999asymptotic}. Such assumptions cannot always
be made, consider for example the cases in
\citep{1555363,prout07}; our result requires much weaker
assumptions and appears to be more general.

More precisely, we consider a family of stochastic processes
$Y^N$, indexed by $N=1,2,...$ over some Polish space $E$; we
assume that the processes have the property that, as $N\to
\infty$ and the initial conditions $y^N(0)\to y(0)$, the
marginals of the process $Y^N(t)$ converge in distribution to
some deterministic $y(t)$, where convergence is at every fixed
time $t$ (see Hypothesis~\ref{ass-1} below). We show that,
under the (mild) assumption that the processes are Feller, this
is sufficient to obtain that any limit point of an invariant
probability of $Y^N$ is an invariant probability of the
deterministic process (and thus its support is included in its
Birkhoff center). Note that we do not need to assume any semi-flow nor continuity
property for the limiting deterministic process.

In the special case where the deterministic process has a
unique limit point $y^*=\limit{t}{\infty}y(t)$ and where the
sequence of invariant probabilities $\Pi^N$ is tight, it
follows immediately that $\Pi^N$ converges to the Dirac mass at
$y^*$. This result is known in the context of stochastic
approximation algorithms; our results here extend it to a more
general setting.

Our result is also more general as it applies to other cases.
Mean field interaction models were often used as practical
approximations of complex interacting object systems, where the
stationary distribution of the system $Y^N$ is approximated by
the stationary regime of an ordinary differential equation
(ODE); this was applied for example to TCP connections
\citep{tinnakornsrisuphap2003limit,baccelli2006http,graham2009interacting},
HTTP flows \citep{baccelli2004metastable}, bandwidth sharing
between streaming and file transfers
\citep{kumar2007integrating}, mobile networks \citep{1555363},
robot swarms \citep{martinoli2002modeling}, transportation
systems \citep{afanassieva1997models}, reputation systems
\citep{le2007generic}, just to name a few. Previous results are
obtained when the ODE has a unique limit point to which all
trajectories converge. Not only does our result extends this
finding to more general spaces, it also extends it to the cases
where there is not a unique limit point. For example, in
\citep{prout07}, the ODE a unique limit point under some
restrictive assumptions on the model parameters; if these
assumptions do not hold, the ODE may have limit cycles, as
shown in \citep{cho-mama-2010}, and nothing can be concluded
from \citep{prout07}. Using our results, it follows that the
limit points of any invariant probability has a support
included in the limit cycles.

\section{Assumptions and Notation}

\subsection{A Collection of Random Processes} Let $(E,d)$ be a Polish space
and $\calP(E)$ the set
of probability measures on $E$, endowed with the topology of
weak convergence. Let $\calC_b(E)$ be the set of bounded
continuous
functions from $E$ to $\Reals$.

We are given a collection of probability spaces
$(\Omega^N,\calF^N, \P^N)$ indexed by $N =1,2,3,...$ and for
every $N$ we have a process $Y^N$ defined on
$(\Omega^N,\calF^N, \P^N)$. Time is either discrete or
continuous. In the discrete time case, $Y^N(t)$ is a collection
of random variables indexed by $N=1,2,3...$ and $t \in \Nats$.
In the continuous time case, let $D_E[0,\infty)$ be the set of
c\'adl\'ag functions $[0,\infty) \to E$; $Y^N$ is then a
stochastic process with sample paths in $D_{E}[0,\infty)$.

We denote by $Y^N(t)$ the random value of $Y^N$ at time
$t\geq 0$. Let $E^N\subset E$ be the support of $Y^N(0)$, so
that $\P^N(Y^N(0)\in E^N)=1$.

We assume that, for every $N$, the process $Y^N$ is Feller, in
the sense that for every $t\geq 0$ and $h \in \calC_b(E)$,
$\E^N\left. \lb h(Y^N(t))\right| Y^N(0)=y_0\rb$ is a
continuous function of $y_0 \in E$.

\begin{defi} A probability $\Pi^N\in \calP(E)$ is invariant for $Y^N$ if
$\Pi^N(E^N)=1$ and for every $h \in \calC_b(E)$ and every $t
\geq 0$:
  \ben
 \int_E \E^N\left. \lb h\lp Y^N(t)\rp\right|Y^N(0)=y\rb \Pi^N(dy) = \int_E h(y)
 \Pi^N(dy)
  \een\label{def1}
  \end{defi}

\subsection{A Deterministic Measurable Process}

Further, let $\varphi$ be a deterministic process, i.e. a
mapping
 \ben \barr{rccl}
 \varphi:& \T \times E &\to &E
 \\
&t,y_0 & \mapsto & \varphi_t(y_0)
 \earr\een where $\T=\Nats$ or $\T=[0,\infty)$.

We assume that $\varphi_t$ is measurable for every fixed time
$t\geq 0$. Note that there is no assumption here that $\varphi$ is continuous nor that it is a flow.

\begin{defi}\label{def2}
A probability $\Pi\in \calP(E)$ is invariant for $\varphi$ if for
every $h \in \calC_b(E)$ and every $t \geq 0$:
  \ben
 \int_E h\lp \varphi_t(y)\rp \Pi(dy) = \int_E h(y)
 \Pi(dy)
  \een
  \end{defi}

%

\subsection{Convergence Hypothesis} We assume that, for every fixed $t$ the processes
$Y^N$ converge in distribution to the deterministic process
$\varphi$ as $N \to \infty$ for every collection of converging
initial conditions. More precisely:

\begin{assh}\label{ass-1}For every $y_0$ in $E$, every sequence
$(y^N_0)_{N=1,2,...}$ such that $y^N_0  \in E^N$ and
$\limit{N}{\infty}y^N_0=y_0$, and every $t\geq0$, the
conditional law of $Y^N(t)$ given $Y^N(0)=y^N_0$ converges
weakly to the Dirac mass at $\varphi_t(y_0)$. That is
$$\limit{N}{\infty} \E^N\left. \lb h(Y^N(t))\right| Y^N(0)=y_0^N\rb = h \circ \varphi_t(y_0)$$ for all
$h \in \calC_b(E)$ and any fixed $t\geq 0$.
\end{assh}

\subsection{Examples.} In discrete time, $Y^N$ is a Markov chain
on $E^N$, as in \citep{le2007generic}, where $E=\calP(S)$ for
some compact set $S$, $Y^N$ is the occupancy measure of a
process on $S$, and $E^N$ is the set of probabilities that are
the sum of $N$ Dirac masses. Here Definition~\ref{def1}
coincides with invariant probability for a Markov chain. The
deterministic process is an iterated map, and
Definition~\ref{def2} coincides with invariant probability of
an iterated map.

In continuous time, $Y^N$ may be a Markov process on $E^N$, as
in
\citep{kurtz1970solutions,sandholm2006population,prout07,graham1994chaos,
benaim2008cmf}. Definition~\ref{def1} coincides here with
invariant probability for a Markov process. The deterministic
process is a semi-flow, and Definition~\ref{def2} coincides
with invariant probability for semi-flows. If $E$ is finite
dimensional, the deterministic process is an ODE or a
differential inclusion.

Still in continuous time, $Y^N$  may also be the continuous
linear interpolation of a discrete time process, as in
\citep{prout07,benaim2008cmf} (in this case it is not a Markov
process). An invariant probability for $Y^N$ is here an
invariant probability of the interpolated Markov chain.

Hypothesis~\ref{ass-1} holds in
\citep{le2007generic,kurtz1970solutions,sandholm2006population,prout07,
benaim2008cmf} as a consequence of stronger convergence
results; for example in \citep{kurtz1970solutions} there is
almost sure, uniform convergence for all $t \in [0,T]$, for any
$T\geq 0$.

\section{Convergence of Invariant Probabilities}

\subsection{Main Theorem}
\begin{thm}Assume Hypothesis~\ref{ass-1} holds and let $\Pi \in \calP(E)$ be a limit point of the
sequence $\Pi^N$, where $\Pi^N$ is an invariant probability for
$Y^N$. Then $\Pi$ is an invariant probability for $\varphi$.
\end{thm}
\begin{mypr}
Let $N_k$ be a subsequence such that
$\limit{k}{\infty}\Pi^{N_k}=\Pi$ in the weak topology on
$\calP(E)$. By Skorohod's representation theorem for Polish
spaces \citep[Thm 1.8]{ethier-kurtz-05}, there exists a common
probability space $(\Omega, \calF,\P)$ on which some random
variables $X^k$ for $k \in \Nats$ and $X$ are defined such that
 \ben
 \bracket{
 \mbox{ law of } X^{k} = \Pi^{N_k}\\
 \mbox{ law of } X = \Pi\\
 X^{k}\to X \; \P-\mbox{a.s.}
  }
 \een

Fix some $t \geq 0$ and $h \in \calC_b(E)$, and define, for $k
\in \Nats$ and $y \in E$
  \ben
  a^k(y) \eqdef \espc{h\lp Y^{N_k}(t)\rp}{Y^{N_k}(0)=y}
  \een

Since $\Pi^{N_k}$ is invariant for $Y^{N_k}$:
  \be
  \int_E a^k(y) \Pi^{N_k}(dy) = \int_E h(y)
  \Pi^{N_k}(dy)\label{eq-877}
  \ee

Hypothesis~\ref{ass-1} implies that
$\limit{k}{\infty}a^k(x^k)=h(\varphi_t(x))$ for every sequence
$x^k$ such that $x^k\in E^{N_k}$ and $\limit{k}{\infty}x^k=x
\in E$. Now $X^k \in E^{N_k}$ $\P-$ almost surely, since the law
of $X^k$ is $\Pi^{N_k}$ and $\Pi^{N_k}$ is invariant for
$Y^{N_k}$. Further, $X^k\to X$ $\P-$ almost surely; thus

 \be
  \limit{k}{\infty} a^k(X^k) = h(\varphi_t(X)) \;\;\;
  \P-\mbox{ almost surely}
  \ee
 Now $a^k(X^k)\leq \norm{h}_{\infty}$ and, thus, by dominated
 convergence:

 \be
 \limit{k}{\infty}\esp{a^k(X^k)}=\esp{h(\varphi_t(X))}
 \ee
 Using \eref{eq-877}:
  \be
 \limit{k}{\infty}\int_E h(y)
  \Pi^{N_k}(dy) = \int_E h\lp\varphi_t(y)\rp
  \Pi(dy)
  \ee and thus
\be \int_E h(y) \Pi(dy)= \int_E h\lp\varphi_t(y)\rp
  \Pi(dy)\ee This holds for any $h\in \calC_b(E)$ and
$t\geq 0$, which shows that $\Pi$ is invariant for $\varphi$.
 \end{mypr}

\subsection{The Continuous Semi-Flow Case} Note that our assumptions on $\varphi$ are very
weak. We now make an additional assumption:
 \begin{defi}\label{def3}
The deterministic process $\varphi$ is a continuous semi-flow
if \begin{enumerate}
\item $\varphi_0(y) = y$
     \item $\varphi_{s+t}=\varphi_s\circ
\varphi_t$ for all nonnegative $s$ and $t$
     \item  $\varphi_{t}(y)$ is continuous in $t$ and  $y$
   \end{enumerate}
 \end{defi}

If $\varphi$ is a continuous semi-flow, it follows from
 Poincar\'e's recurrence theorem \citep{mane} that the support
 of any limit point of $\Pi^N$ is included in the closure of the recurrent set :
 $$R(\varphi) = \{x \in E: \: \liminf_{t \rightarrow \infty} d(x,\varphi_t(x)) = 0\}.$$

 In particular, if the semi-flow has a
 unique limit point, we have:

 \begin{coro}\label{coro-uni}Assume Hypothesis~\ref{ass-1} holds and that $\varphi$ is a continuous semi-flow.
 Let
$\Pi^N$ be a sequence of invariant probabilities for $Y^N$.
Assume that
  \begin{enumerate}
    \item the sequence $(\Pi^N)_{N=1,2,...}$ is tight;
    \item there is some $y^*\in E$ such that for all $y \in
        E$, $\limit{t}{\infty}\varphi_t(y)=y^*$.
  \end{enumerate}
  It follows that the sequence $\Pi^N$ converges weakly to the Dirac mass at
  $y^*$.
 \end{coro}

 Recall that tightness means that for every $\eps>0$ there is
 some compact set $K \subset E$ such that $\Pi^N(K)\geq 1-\eps$
 for all $N$. If $E$ is compact then $(\Pi^N)_{N=1,2,...}$ is
 automatically tight.%
%
%

\bibliographystyle{elsarticle-harv}
\bibliography{leb}

\end{document}